\label{key}
\documentclass[12pt, reqno]{amsart}
\usepackage{amssymb, amsthm, amsmath, amsfonts}
\usepackage{array, epsfig}
\usepackage{bbm}
\usepackage{hyperref}
\usepackage[numbers,square]{natbib}
\usepackage{color}
\usepackage{graphicx}
\usepackage{multirow}

  \evensidemargin
\oddsidemargin
\newtheorem{theorem}{Theorem}[section]
\newtheorem{corollary}[theorem]{Corollary}
\newtheorem{remark}{Remark}
\newtheorem*{example}{Example}
\newtheorem{lemma}[theorem]{Lemma}

\newcommand{\alg}{\mathbb{A}}
\newcommand{\R}{\mathbb{R}}
\newcommand{\Z}{\mathbb{Z}}
\newcommand{\C}{\mathbb{C}}
\newcommand{\N}{\mathbb{N}}
\newcommand{\Q}{\mathbb{Q}}

\newcommand{\bx}{\mathbf{x}}
\newcommand{\ba}{\mathbf{a}}

\newcommand{\dd}{{\rm d}}

\newcommand{\ind}{\mathbbm{1}}
\newcommand{\Sl}{Sal}
\newcommand{\Pl}{\mathcal{P}}
\newcommand{\Per}{Per}

\newcommand{\T}{\mathbb{T}}

\newcommand{\h}{Q} 

\newcommand{\pderiv}[2]{\frac{\partial #1}{\partial #2}}
\newcommand{\jac}[3]{J_{#1}^{(#2,#3)}}

\newcommand{\Vol}{\mathop{\mathrm{Vol}}\nolimits}

\newcommand{\Pf}{\operatorname{Pf}}
\newcommand{\sign}{\operatorname{sign}}
\newcommand{\sgn}{\operatorname{sgn}}

\makeatletter
\newcommand*{\house}[1]{%
  \mathord{%
    \mathpalette\@house{#1}%
  }%
}
\newcommand*{\@house}[2]{%
  \dimen@=\fontdimen8 %
      \ifx#1\scriptscriptstyle\scriptscriptfont
      \else\ifx#1\scriptstyle\scriptfont
      \else\textfont\fi\fi
      3 %
  \sbox0{%
    $#1%
      \vrule width\dimen@\relax
      \overline{%
        \kern2\dimen@
        \begingroup 
          #2%
        \endgroup
        \kern2\dimen@
      }%
      \vrule width\dimen@\relax
      \mathsurround=1.5\dimen@ 
    $%
  }%
  \ht0=\dimexpr\ht0-\dimen@\relax
  \dp0=\dimexpr\dp0+2\dimen@\relax
  \vbox{%
    \kern\dimen@ 
    \copy0 %
  }%
}

\title{On the distribution of Salem numbers}

\thanks{The research of the first author has been supported by project SFB 1283 at Bielefeld University (Germany). The research of the second author has been supported by project IRTG 2235 at Bielefeld University (Germany).}
	
\author{Friedrich G\"otze}
\address{Friedrich G\"otze: Faculty of Mathematics,
	Bielefeld University,
	P. O. Box 10 01 31,
	33501 Bielefeld, Germany}
\email{goetze@math.uni-bielefeld.de}

\author{Anna Gusakova}
\address{Anna Gusakova: Faculty of Mathematics,
	Bielefeld University,
	P. O. Box 10 01 31,
	33501 Bielefeld, Germany}
\email{agusakov@math.uni-bielefeld.de}

\begin{document}

\maketitle

\begin{abstract}
In this paper we study the problem of counting Salem numbers of fixed degree. Given a set of disjoint intervals $I_1,\ldots, I_{k}\subset \left[0;\pi\right]$, $1\leq k\leq m$ let $\Sl_{m,k}(\h,I_1,\ldots,I_{k})$ denote the set of ordered $(k+1)$-tuples $\left(\alpha_0,\ldots,\alpha_{k}\right)$ of conjugate algebraic integers, such that $\alpha_0$ is a Salem numbers of degree $2m+2$ satisfying $\alpha\leq \h$ for some positive real number $\h$ and $\arg\alpha_i\in I_i$. We derive the following asymptotic approximation 
	\[
	\#\Sl_{m,k}(\h,I_1,\ldots,I_{k})=\omega_m\,\h^{m+1}\,\int\limits_{I_1}\ldots\int\limits_{I_{k}}\rho_{m,k}(\boldsymbol\theta)\dd \boldsymbol\theta+O\left(\h^{m}\right),\quad \h\rightarrow\infty,
	\]
	providing explicit expressions for the constant $\omega_m$ and the function $\rho_{m,k}(\boldsymbol\theta)$.
	 Moreover we derive a similar   asymptotic formula
for the set of all Salem numbers of fixed degree and absolute value bounded by $\h$ as $\h\rightarrow\infty$.
\end{abstract}

\section{Introduction}

Problems concerning  the distribution of algebraic numbers have a long history \cite{GKZ17, Bar14, GG17, MV08, Sch93, vB99}. Recall that an algebraic number $\alpha$ is a complex number such that there exists an irreducible  polynomial $P$ over $\Q$ with integer co-prime coefficients and positive leading coefficient such that $P(\alpha)=0$. This polynomial is called minimal polynomial of the algebraic number $\alpha$ and other roots of $P$ are called Galois conjugates of $\alpha$. Moreover if the leading coefficient of the minimal polynomial equals $1$, $\alpha$ is called algebraic integer.

For the distribution of algebraic numbers we consider sets $\alg_n$ of algebraic numbers with fixed degree $n\in\N$. Usually this set will be countable and dense in $\R$. Hence we shall restrict the counting problem to finite subsets of $\alg_n$ depending on a real parameter $\h > 1$ and ask how the cardinality of this set changes as $\h\rightarrow \infty$.
	As one of the many choices  consider for example the set of algebraic numbers in $\alg_n$ with absolute multiplicative Weil height bounded by $\h$ \cite{Sch93, MV08}, or the set of algebraic numbers in $\alg_n$ with na\"ive height bounded by $\h$ and lying in some fixed set $D\subset\C$ \cite{GKZ15, dK14}. 

In this paper we will consider similar questions for a special subset of algebraic integers, namely the so-called Salem numbers. Let us start with some definitions. A {\it Salem number} is a real algebraic integer $\alpha>1$ such that all its Galois conjugates have absolute value less or equal to $1$ and at least one of them has absolute value equal to $1$. Let $\alpha'$ denote a Galois conjugate of the Salem number $\alpha$ lying on the complex unit circle $\T$, e.g. $\left|\alpha'\right|=1$. Since $\alpha'$ and its complex conjugate $\overline{\alpha'}=\left(\alpha'\right)^{-1}$ are Galois conjugates we conclude that the minimal polynomial $P_{\alpha}$ of a Salem number $\alpha$ is self-reciprocal. Recall that a polynomial $P\in\Z[t]$ of degree $n$ is called self-reciprocal if
\[
P(t)=t^nP\left(t^{-1}\right).
\]
Moreover the polynomial $P_{\alpha}$ is of even degree $2(m+1)$, otherwise $P_{\alpha}(-1)=-P_{\alpha}(-1)=0$ which contradicts to the irreducibility of $P_{\alpha}$. Thus, all Galois conjugates of Salem number $\alpha$ (except for $\alpha^{-1}$) have absolute value $1$ and lie on the unit circle $\T$ in the complex plane. We will denote them by $\alpha_1,\bar{\alpha}_1,\ldots, \alpha_{m},\bar\alpha_{m}\in \T$. We shall use these two properties as a description of the set of Salem numbers. 

Denote by $\Sl_m$ the set of all Salem numbers of degree $2(m+1)$. Our aim is to describe the distribution of Salem numbers by considering some finite subsets of $\Sl_m$ with given properties and investigating how the cardinality of those sets depend on the chosen parameters.

It should be mentioned that Salem numbers play an important role in many areas of mathematics, such as number theory, algebra and dynamical systems. For more details we refer to the papers \cite{PiSa, Boy81, Sal63, EvWar}. In particular the smallest Salem number is closely related to  Lehmer's conjecture \cite{Leh}.

\subsection{Counting Salem numbers}

Given some real $\h>1$ and integer $m$ let us introduce the following finite subset
\[
\Sl_m(\h):=\left\{\alpha\in\Sl_m\colon \alpha\leq\h\right\},
\]
which consists of all Salem numbers of degree $2(m+1)$ lying in the interval $(1;\h]$. In this subsection we introduce an asymptotic formula for the cardinality of the set $\Sl_m(\h)$ as $\h\rightarrow \infty$.

\begin{theorem}\label{th1}
	For any integer $m$ we have
	\[
	\#\Sl_m(\h)=\omega_m\,\h^{m+1}+O\left(\h^{m}\right),\quad \h\rightarrow\infty,
	\]
	where
	\begin{equation}\label{eq_70}
	\omega_m:=\frac{2^{m(m+1)}}{m+1}\,\prod\limits_{k=0}^{m-1}\frac{k!^2}{(2k+1)!}.
	\end{equation}
\end{theorem}

It should be mentioned that the set of Salem numbers is a subset of a more general class of algebraic integers, namely Perron numbers. A Perron number is a real algebraic integer $\alpha>1$ such that all its Galois conjugates have absolute value less then $\alpha$. The problem of counting Perron numbers has been studied by F. Calegari and Z. Huang \cite{CH17}. 
 Denote by $\Per_n$ the set of all Perron numbers of degree $n$ and for some real $\h>1$ define the following finite subset
 \[
 \Per_n(\h):=\left\{\alpha\in\Per_n\colon \alpha\leq\h\right\}.
 \]
 Then for any natural $n$ the following asymptotic approximation holds
 \[
 \# \Per_n(\h)= d_n\h^{n(n+1)/2}+O(\h^{n(n-1)/2}),\quad \h\rightarrow\infty,
 \]
where
 \[
 d_{2s}=\frac{1}{2s+1}\prod\limits_{k=0}^{s-1}\left(\frac{k!^22^{2k+1}}{(2k+1)!}\right)^2;\quad  d_{2s+1}=\frac{2^{2s+1}s!^2}{(2s+1)!}d_{2s}.
 \]

The proof methods used in \cite{CH17} do not apply to the case of Salem numbers, but we use some modification of those arguments to prove Theorem \ref{th1}.

\subsection{Salem numbers with given distribution of their Galois conjugates}

In this subsection we consider a slightly more general problem. 

Given some real $\h>1$, integer $m$ and disjoint intervals $I_1,\ldots, I_{k}\subset \left[0;\pi\right]$, $1\leq k\leq m$ denote by $\Sl_{m,k}(\h,I_1,\ldots,I_{k})$ the set of ordered $(k+1)$-tuples $\left(\alpha_0,\ldots,\alpha_{k}\right)\in \R\times\T^{k}$ of conjugate algebraic integers, such that $\alpha_0\in\Sl_m(\h)$ and $\arg\alpha_i\in I_i$ for $1\leq i\leq k$. As in the previous subsection we try to determine the cardinality of this set. Theorem \ref{th2} below provides an asymptotic formula for $\#\Sl_{m,k}(\h,I_1,\ldots,I_{k})$ as $\h$ tends to infinity.

Before we introduce our main result let us consider some additional notations and definitions. Let $A=(a_{i,j})_{i,j=1,\ldots,2n}$ be a skew-symmetric $2n\times 2n$ matrix, which means that $A^{\top}=-A$. The Pfaffian of $A$ is defined by
\[
\Pf(A)=\frac{1}{2^nn!}\sum\limits_{\sigma\in S_{2n}}\sgn(\sigma)\prod\limits_{i=1}^{n}a_{\sigma(2i-1),\sigma(2i)},
\]
where $S_{2n}$ is the symmetric group of the dimension $(2n)!$ and $\sgn(\sigma)$ is the signature of $\sigma$.Recall the following useful formula connecting Pfaffian and determinant of the matrix $A$
\begin{equation}\label{eq_13}
\Pf(A)^2=\det(A).
\end{equation}
Furthermore, let us  introduce the family of classical orthogonal polynomials, called Jacobi polynomials, via
\begin{equation}\label{eq_51}
\jac{n}{a}{b}(t)=\frac{\Gamma\left(a+n+1\right)}{n!\Gamma\left(a+b+n+1\right)}\sum\limits_{j=0}^{n}\binom{n}{j}\frac{\Gamma\left(a+b+n+j+1\right)}{\Gamma\left(a+j+1\right)}\left(\frac{t-1}{2}\right)^j.
\end{equation}
These polynomials are orthogonal to each other with respect to weight function $(1-t)^{a}(1+t)^{b}$ on the interval $\left[-1;1\right]$. For more details we refer reader to \cite{Szego} and Appendix \ref{jacobi}.

\begin{theorem}\label{th2}
	For any integer $m$ and any disjoint intervals $I_1,\ldots, I_{k}\subset \left[0;\pi\right]$, $1\leq k\leq m$ we have
	\[
	\#\Sl_{m,k}(\h,I_1,\ldots,I_{k})=\omega_m\,\h^{m+1}\,\int\limits_{I_1}\ldots\int\limits_{I_{k}}\rho_{m,k}(\boldsymbol\theta)\dd \boldsymbol\theta+O\left(\h^{m}\right),\quad \h\rightarrow\infty,
	\]
	where $\omega_m$ is defined by \eqref{eq_70}. Moreover, the function $\rho_{N,k}(\boldsymbol\theta)$ can be written in the following form
	\[
	\rho_{N,k}(\theta_1,\ldots,\theta_{k})=\prod\limits_{l=1}^{k}\sin\theta_l\Pf\Big[K_N(\cos\theta_i,\cos\theta_j)\Big]_{i,j=1,\ldots,k},
	\]
	where 
	\begin{equation}\label{eq_23}
	K_N(x,y):= \left(\begin{matrix}
	I_N(x,y) & S_N(y,x)\\
	-S_N(x,y) & -D_N(x,y)
	\end{matrix}
	\right)
	\end{equation}
	and for $c:=\left(N\mod{2}\right)$ we write
	\begin{align}
	S_N(x,y):&=\sum\limits_{j=0}^{(N-c)/2-1}\frac{1}{r_j}\left(\psi'_{2j+1}(x)\psi_{2j}(y)-\psi'_{2j}(x)\psi_{2j+1}(y)\right)+\frac{c(N+1)}{4}\jac{N-1}{1}{1}(x),\label{eq_53}\\
	D_N(x,y):&=-\pderiv{S_N(x,y)}{y},\label{eq_54}\\
	I_N(x,y):&=\frac12\int\limits_{-1}^{1}\sign(x-\xi)S_N(\xi,y)\dd\xi-\frac12\sign(x-y)-\frac{c}{2}\jac{N}{0}{0}(x),\label{eq_55}
	\end{align}
where
	\begin{align}
	\psi_{2j}(t)&:=\frac{2}{2j+2+c}\left(\jac{2j+1+c}{0}{0}(t)-c\right),\label{eq_56}\\
	\psi_{2j+1}(t)&:=(t^2-1)\jac{2j+c}{1}{1}(t),\label{eq_57}\\
	r_j&:=\frac{8(2j+1+c)}{(4j+3+2c)(2j+2+c)}.\label{eq_58}
	\end{align}
\end{theorem}

\begin{remark}
	According to the definition of the Pfaffian and the kernel $K_N(x,y)$ it is easy to see that $\rho_{N,k}(\theta_1,\ldots,\theta_{k})$ is a polynomial in $\sin\theta_i$ and $\cos\theta_i$, $1\leq i\leq k$. Moreover, it can be written as
	\[
	\rho_{N,k}(\theta_1,\ldots,\theta_{k})=\prod\limits_{l=1}^{k}\sin\theta_l\prod\limits_{1\leq i<j\leq k}\left|\cos\theta_i-\cos\theta_j\right|g_{N}(\cos\theta_1,\ldots,\cos\theta_{k}),
	\]
	where $g_{N}(x_1,\ldots,x_k)$ is a polynomial.
\end{remark}

\begin{remark}
	It should be also noted that the function $K_N(x,y)$ defined by \eqref{eq_23} - \eqref{eq_58} coincides with the Kernel function of some random matrix ensemble, namely the Jacobi $\beta$-ensemble with $\beta=1$ (see Appendix \ref{theory} and \ref{jacobi}).
\end{remark}

In general, the formula for $\rho_{m,k}(\boldsymbol\theta)$ seems quite complicated, but it may be simplified for $k=1$ and $k=m$.

\begin{corollary}\label{cr3}
	For any integers $m$ we have
	\[
	\rho_{m,1}(\theta)=\sin\theta\,S_{m}(-\cos\theta, -\cos\theta),
	\]
	and
	\[
	\rho_{m, m}(\theta_1,\ldots,\theta_{m})=\prod\limits_{l=1}^{m}\sin\theta_l\prod\limits_{1\leq i<j\leq m}\left|\cos\theta_i-\cos\theta_j\right|.
	\]
\end{corollary}

Corollary \ref{cr3} immediately follows from the proof of Theorem \ref{th2} and equation \eqref{eq_13}.

\begin{figure}[h]
\begin{minipage}[h]{0.49\linewidth}
\center{\includegraphics[width=0.8\linewidth]{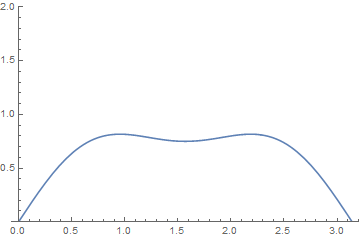} \\$m=2$}
\end{minipage}
\hfill
\begin{minipage}[h]{0.49\linewidth}
\center{\includegraphics[width=0.8\linewidth]{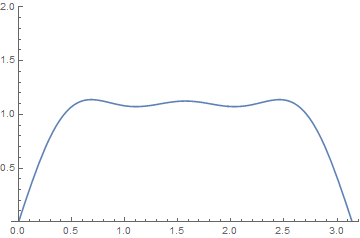} \\$m=3$}
\end{minipage}
\hfill
\hfill
\begin{minipage}[h]{0.49\linewidth}
\center{\includegraphics[width=0.8\linewidth]{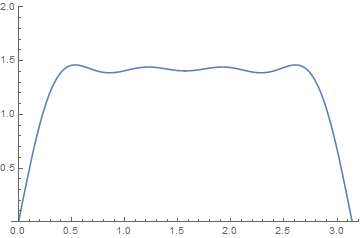} \\$m=4$}
\end{minipage}
\hfill
\begin{minipage}[h]{0.49\linewidth}
\center{\includegraphics[width=0.8\linewidth]{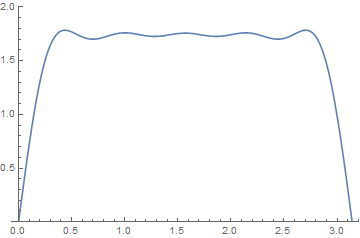} \\$m=5$}
\end{minipage}
\caption{A plot of  $\rho_{m,1}(\theta)$.}
\label{image1}
\end{figure}

\begin{example}
As an example the following density functions may be written explicitly as 
\begin{align*}
\rho_{2,1}(\theta)&=\frac34\sin\theta\left(\cos^2\theta+1\right),\\
\rho_{3,1}(\theta)&=\frac{3}{8}\sin\theta\left(5\cos^4\theta+3\right),\\
\rho_{4,1}(\theta)&=\frac{5}{32}\sin\theta\left(35\cos^6\theta-21\cos^4\theta+9\cos^2\theta+9\right),
\end{align*}
with corresponding plots given in Figure \ref{image1}.
\end{example}

The structure of the paper is the following. In Section 2 we present the proof of Theorem \ref{th1} and Theorem \ref{th2}. Section 3 is devoted to the proof of auxiliary lemmas. In Appendix \ref{theory} and Appendix \ref{jacobi} we collected some facts about the distribution of the eigenvalues of random matrix ensembles needed for the proof of Theorem \ref{th2}.

\section{Proof of Theorem \ref{th1} and Theorem \ref{th2}}

In this section we present the proofs of Theorem \ref{th1} and Theorem \ref{th2} simultaneously since they differ in some details only.

Denote by $\Pl_m(\h)$ the set of self-reciprocal monic polynomials $P\in\Z[t]$ of degree $2(m+1)$ having $2m$ roots lying on the unit circle $\T$ and two positive real roots $\alpha$, $\alpha^{-1}$ bounded by $\h$. Moreover let $\Pl^{I}_m(\h)$ denote the subclass of irreducible polynomials $P\in\Pl_m(\h)$ and let  $\Pl^{R}_m(\h)$ denote the subclass of reducible polynomials $P\in\Pl_m(\h)$. By definition $\Pl^{I}_m(\h)$ is the set of minimal polynomials of Salem numbers $\alpha\in\Sl_m(\h)$ and thus
\begin{equation}\label{eq_12}
\#\Sl_m(\h) = \#\Pl^{I}_m(\h).
\end{equation}

Given a polynomial $P\in\R[t]$ denote by $\mu_P$ a counting measure for the roots of $P$ lying on the unit circle
\begin{equation}\label{eq_19}
\mu_P:=\sum\limits_{\theta\colon P\left(e^{i\theta}\right)=0}\delta_\theta,
\end{equation}
where $\delta_\theta$ is the unit point mass in $\theta$.
Then, for any set of disjoint intervals $I_1,\ldots, I_k\subset[0;\pi]$ and a polynomial $P\in \R[t]$ the following quantity
\begin{equation}\label{eq_20}
\mu_P(I_1,\ldots,I_k):=\prod\limits_{i=1}^k\mu_P(I_i)
\end{equation}
is the number of ordered $k$-tuples $(\theta_1,\ldots,\theta_k)\in I_1\times\ldots\times I_k$ such that 
\[
P(e^{i\theta_1})=\ldots=P(e^{i\theta_k})=0.
\]
Then, obviously,
\begin{equation}\label{eq_16}
\#\Sl_{m,k}(\h,I_1,\ldots,I_{k})=\sum\limits_{l=0}^{\infty}l\cdot\#\left\{P\in \Pl^{I}_m(\h)\colon \mu_P(I_1,\ldots,I_k)=l\right\}
\end{equation}
and our problem reduces to counting integral irreducible polynomials with a prescribed root distribution.

Our approach in this case will be to consider the set of polynomials with real instead of integer coefficients. Thus, identifying a polynomial of degree $d$ with the vector of its coefficients as a point in $\R^{d+1}$ we transform the algebraic problem to the geometric problem of counting a number of integer points (points with integer coefficients) inside specific sets in $\R^{d+1}$.

Following this idea denote by $V_m\subset \R^{m+1}$ a set of points $\ba:=(a_1,\ldots, a_{m+1})$ such that the roots of polynomial 
\[
P_{\ba}(t):=t^{2m+2}+a_1t^{2m+1}+\ldots+a_{m+1}t^{m+1}+\ldots+a_1t+1
\]
have the following form 
\begin{align}
e^{i\theta_1},e^{-i\theta_1},\ldots,e^{i\theta_{m}}, e^{-i\theta_{m}}&\in\T,\notag\\
y&\in (1;\h],\quad y^{-1}, \label{eq_1}
\end{align}
for some $0\leq\theta_1\leq\ldots\leq\theta_{m}\leq\pi$. Moreover, denote by $V^l_{m,k}\subset V^m$ the set of points $\ba\in V_m$ such that the polynomial $P_{\ba}(t)$ satisfies the additional condition
\[
\mu_{P_{\ba}}(I_1,\ldots, I_k)=l.
\]
In order to simplify notation in the rest of the proof we will write $\Sl_{m,k}(\h)$ instead of $\Sl_{m,k}(\h,I_1,\ldots,I_{k})$ and will assume $V_m=V^1_{m,0}$.

Denote by $v^l_{m,k}$ the volume of the set $V^l_{m,k}$ given by
\begin{equation}\label{eq_4}
v^l_{m,k}=\int\limits_{\R^{m+1}}\ind_{V^l_{m,k}}(\bx)\dd\bx,
\end{equation}
where $\ind_{B}(\cdot)$ denotes the indicator function of a set $B\subset \R^{m+1}$.

According to our convention every polynomial from $\Pl_m(\h)$ represents an integer point in $V_m$. Hence, the first step is to  count the integer points in the set $V^l_{m,k}$. We will show that for any integer $0\leq k\leq m$ and $l\ge 1$ this number is asymptotically equal to the volume $v^l_{m,k}$ of $V^l_{m,k}$ as $\h\rightarrow \infty$.

Approximating the number of integer points in a large set $B\subset\R^{m+1}$ by the volume  of this set is a classical approach. One of the earliest references we are aware of in this direction are due to Lip\-schitz~\cite{rL65} and Davenport \cite{hD51}. In order to get a good estimate one needs to impose some regularity conditions on the boundary of $B$. In accordance with~\cite[Definition 2.2]{mW12}, we say that the boundary $\partial B$ of a set $B$ is of Lipschitz class $(M,L)$  if there exist $M$ maps
$\phi_1, \dots, \phi_M : [0; 1]^{m} \to\R^{m+1}$ satisfying a Lipschitz condition
\[
|\phi_i(x) - \phi_i(y)| \le L|x - y| \text{  for  } x, y \in [0;1]^{m}, \quad  i = 1, \dots,M,
\]
such that $\partial B$ is covered by the images of the maps $\phi_i$. 

\begin{lemma}\label{lm3}
	For any integer $0\leq k\leq m$ and $l\ge 1$ the boundary of $V^l_{m,k}$ is of Lipschitz class $(M,c\,\h)$ for some fixed $M$ and $c$ independent of $\h$. Moreover 
	\[
	\left|\#\left(V^l_{m,k}\cap \Z^{m+1}\right)-v^l_{m,k}\right|\leq C_1\h^{m},
	\]
	where $C_1$ depends on $m$, $l$ and $k$ only.
\end{lemma}

For the proof see Subsection \ref{proof_lm3}.

Note that Lemma \ref{lm3} allows us to estimate the number of all integer polynomials in $\Pl_m(\h)$  satisfying the same conditions to the polynomials with real coefficients forming the set $V^l_{m,k}$. But for our purpose we need to count the number of irreducible polynomials only. Thus, our next step is to show that the number of reducible polynomials in $\Pl_m(\h)$ is relatively small and can be estimated by $O(\h^{m})$. The following lemma gives the asymptotic behavior of $\#\Pl^{R}_m(\h)$ as $\h\rightarrow\infty$.

\begin{lemma}\label{lm4}
For some $C_2>0$ depending on $m$ only we have
\[
\#\Pl^{R}_m(\h)\leq C_2\,\h^{m}.
\]
\end{lemma}

The proof of this lemma is given in Subsection \ref{proof_lm4}.

Finally, using Lemma \ref{lm3} and Lemma \ref{lm4} we get
\[
\left|\#\left\{P\in \Pl^{I}_m(\h)\colon \mu_P(I_1,\ldots,I_k)=l\right\}-v^l_{m,k}\right|\leq C_3\,\h^{m},
\]
and from \eqref{eq_12} and \eqref{eq_16} with $v_m:=v_{m,0}^1$ we conclude
\begin{align}
\#\Sl_{m}(\h)&=v_m+O(\h^{m}),\label{eq_17}\\
\#\Sl_{m,k}(\h)&=\sum\limits_{l=0}^{\infty}lv^l_{m,k}+O(\h^{m}).\label{eq_18}
\end{align}

The last step of the proof is devoted to evaluation of $v_m$ and $\sum\limits_{l=0}^{\infty}lv^l_{m,k}$. Let us introduce the following representation of the points $\ba\in V_m$ in terms of the roots \eqref{eq_1} of a polynomial $P_{\ba}(t)$.

\begin{lemma}\label{lm1}
Given a polynomial 
\[
P_{\ba}(t):=t^{2m+2}+a_1t^{2m+1}+\ldots+a_{m+1}t^{m+1}+\ldots+a_1t+1
\]
with roots \eqref{eq_1} we have
\begin{align}
a_{2i-1}&= -\sum\limits_{j=0}^{i-1}{m-2j \choose i-j-1 }\sum\limits_{0\leq k_1<\ldots<k_{2j+1}\leq m}\prod_{s=1}^{2j+1}z_{k_{s}},\label{eq_3}\\
a_{2i}&=\sum\limits_{j=0}^{i}{m-2j+1 \choose i-j}\sum\limits_{0\leq k_1<\ldots<k_{2j}\leq m}\prod_{s=1}^{2j}z_{k_{s}},\label{eq_2}
\end{align}
where $1\leq i\leq \left[\frac{m}{2}\right]+1$, $z_0:=z_0(y)=y+y^{-1}$ and $z_k:=z_k(\theta_k)=2\cos\theta_k$, $1\leq k\leq m$.
\end{lemma}

The proof of this lemma is given in Subsection \ref{proof_lm1}.

Let us consider the following simplex
\begin{equation}\label{eq_26}
L_{m}:=\left\{\bx\in \R^{m}\colon 0\leq x_1\leq x_2\leq\ldots\leq x_{m}\leq \pi \right\}.
\end{equation}
In Lemma \ref{lm1} we defined a bijective map $f: (1,\h] \times L_{m}\rightarrow V_m$ via \eqref{eq_3} and \eqref{eq_2} as
\[
f_i(y,\theta_1,\ldots,\theta_{m}):=a_i.
\]
Using this mapping we can rewrite \eqref{eq_4} as follows
\[
v_m=\int\limits_{1}^{\h}\int\limits_{L_{m}}J_f(y,\theta_1,\ldots,\theta_{m})\dd y\,\dd\theta_1\ldots\dd\theta_{m},
\]
where 
\begin{equation}\label{eq_5}
J_f(y,\theta_1,\ldots,\theta_{m})=\left|\frac{\partial\left(f_1,\ldots,f_{m+1}\right)}{\partial\left(y,\theta_1,\ldots,\theta_{m}\right)}\right|
\end{equation}
is the Jacobian of the map $f$.

\begin{lemma}\label{lm2}
Given the map $f: (1,\h]\times L_{m}\rightarrow V_m$ defined by \eqref{eq_3} and \eqref{eq_2} we have
\[
J_f(y,\theta_1,\ldots,\theta_{m})=2^{\frac{m(m+1)}{2}}\left(1-\frac{1}{y}\right)\prod\limits_{l=1}^{m}\left(y+\frac{1}{y}-2\cos\theta_l\right)\prod\limits_{l=1}^{m}\sin\theta_l\prod\limits_{1\leq i<j\leq m}\left|\cos\theta_i-\cos\theta_j\right|.
\]
\end{lemma}

The evaluation of the Jacobian is postponed to Subsection \ref{proof_lm2}.

Since  $J_f(y,\theta_1,\ldots,\theta_{m})$ is invariant with respect to permutations of $\theta_1,\ldots,\theta_{m}$, we may write
\[
v_m=\frac{1}{m!}\int\limits_{1}^{\h}\int\limits_{0}^{\pi}\ldots\int\limits_{0}^{\pi}J_f(y,\theta_1,\ldots,\theta_{m})\dd y\,\dd\theta_1\ldots\dd\theta_{m}.
\]
Now using Lemma \ref{lm2} together with the change of variables 
\begin{align*}
x_0&=y+y^{-1},\\
x_i&=-\cos\theta_i,\quad 1\leq i\leq m,
\end{align*}
\begin{align}
v_m&=\frac{2^{m(m+1)/2}}{m!}\int\limits_{2}^{\h+\h^{-1}}\int\limits_{-1}^{1}\ldots\int\limits_{-1}^{1}\prod\limits_{l=1}^{m}\left(x_0+2x_l\right)\prod\limits_{1\leq i<j\leq m}\left|x_i-x_j\right|\dd x_0\,\dd x_1\ldots\dd x_{m} \notag \\
&=\omega_{m}\left(\h+\h^{-1}\right)^{m+1}+\omega_{m-1}\left(\h+\h^{-1}\right)^{m}+\ldots+\omega_0\left(\h+\h^{-1}\right)-\sum\limits_{i=0}^{m}2^{i+1}\omega_{i},\label{eq_8}
\end{align}
where
\begin{equation}\label{eq_9}
\omega_l:=2^{m-l}\kappa_m\,(l+1)^{-1}\int\limits_{-1}^{1}\ldots\int\limits_{-1}^{1}\sigma_{m-l}\left(x_1,\ldots,x_{m}\right)\prod\limits_{1\leq i<j\leq m}\left|x_i-x_j\right|\dd x_1\ldots\dd x_{m},
\end{equation}
and $\sigma_{k}$ denotes the $k$-th elementary symmetric polynomial.

It is easy to see that the integrals in \eqref{eq_9} are Selberg's type integrals. In particular, taking $y_i=\frac{x_i+1}{2}$ for $l=m$ we obtain 
\begin{align*}
\omega_{m}&=\frac{2^{m(m+1)/2}}{(m+1)!}\,\int\limits_{-1}^{1}\ldots\int\limits_{-1}^{1}\prod\limits_{1\leq i<j\leq m}\left|x_i-x_j\right|\dd x_1\ldots\dd x_{m}\\
&=\frac{2^{m(m+1)}}{(m+1)!}\,\int\limits_{0}^{1}\ldots\int\limits_{0}^{1}\prod\limits_{1\leq i<j\leq m}\left|y_i-y_j\right|\dd y_1\ldots\dd y_{m},
\end{align*}
which is a special case of Selberg's integral formula \cite{Sel44}
\begin{align*}
S_{n}\left(\alpha,\beta,\gamma\right):&=\int\limits_{0}^{1}\ldots\int\limits_{0}^{1}\prod\limits_{i=1}^nt_i^{\alpha-1}(1-t_i)^{\beta-1}\prod\limits_{1\leq i<j\leq n}\left|t_i-t_j\right|^{2\gamma}\dd t_1\ldots\dd t_{n} \notag\\
&=\prod\limits_{j=0}^{n-1}\frac{\Gamma\left(\alpha+j\gamma\right)\Gamma\left(\beta+j\gamma\right)\Gamma\left(1+(j+1)\gamma\right)}{\Gamma\left(\alpha+\beta+(n+j-1)\gamma\right)\Gamma\left(1+\gamma\right)},
\end{align*}
for $\alpha=\beta=1$ and $\gamma = \frac12$. Thus, we conclude
\[
\omega_{m}=\frac{2^{m(m+1)}}{(m+1)!}\,\prod\limits_{j=0}^{m-1}\frac{\Gamma\left(1+\frac12j\right)^2\Gamma\left(\frac32+\frac12j\right)}{\Gamma\left(1+\frac12(m+1+j)\right)\Gamma\left(\frac32\right)}=\frac{2^{m(m+1)}}{m+1}\,\prod\limits_{k=0}^{m-1}\frac{k!^2}{(2k+1)!}.
\]
Substituting this into \eqref{eq_8} leads to
\begin{equation}\label{eq_11}
v_m=\omega_m\,\h^{m+1}+O(\h^{m}),
\end{equation}
which together with \eqref{eq_17} finishes the proof of Theorem \ref{th1}.

The evaluation of 
\[
W:=\sum\limits_{l=0}^{\infty}lv^l_{m,k}
\]
is a bit more involved. Since it does not seem possible to derive compact representations for every $v^l_{m,k}$ separately, we shall look for representation of the whole sum $W$ instead. 

First of all, using \eqref{eq_4} and the definition of the sets $V^l_{m,k}$ we conclude
\[
W=\sum\limits_{l=0}^{\infty}l\cdot \int\limits_{\R^{m+1}}\ind_{V^l_{m,k}}(\bx)\dd\bx=\int\limits_{\R^{m+1}}\sum\limits_{l=0}^{\infty}l\cdot \ind_{V^l_{m,k}}(\bx)\dd\bx= \int\limits_{\R^{m+1}}\mu_{P_{\bx}}(I_1,\ldots,I_k) \ind_{V_{m}}(\bx)\dd\bx.
\]
Applying Lemma \ref{lm1} and using the representation of $\mu_{P_{\bx}}(I_1,\ldots,I_k)$ via \eqref{eq_19} and \eqref{eq_20} we get
\[
W=\int\limits_{1}^{\h}\int\limits_{L_{m}}J_f(y,\theta_1,\ldots,\theta_{m})\prod\limits_{i=1}^k\sum\limits_{j=1}^{m}\delta_{\theta_j}(I_i)\dd y\,\dd\theta_1\ldots\dd\theta_{m},
\]
and since our integrand is invariant with respect to permutations of $\theta_1,\ldots,\theta_{m}$ we write
\[
W=\frac{1}{m!}\int\limits_{1}^{\h}\int\limits_{0}^{\pi}\ldots\int\limits_{0}^{\pi}J_f(y,\theta_1,\ldots,\theta_{m})\prod\limits_{i=1}^k\sum\limits_{j=1}^{m}\delta_{\theta_j}(I_i)\dd y\,\dd\theta_1\ldots\dd\theta_{m}.
\]
From the identity 
\[
\delta_{\theta_i}(I_j)=\ind_{I_j}(\theta_i)
\]
and the fact that intervals $I_1,\ldots, I_k$ are disjoint it follows 
\begin{align*}
W&=\frac{1}{m!}\int\limits_{1}^{\h}\int\limits_{0}^{\pi}\ldots\int\limits_{0}^{\pi}J_f(y,\theta_1,\ldots,\theta_{m})\prod\limits_{i=1}^k\sum\limits_{j=1}^{m}\ind_{I_i}(\theta_j)\dd y\,\dd\theta_1\ldots\dd\theta_{m}\\
&=\frac{1}{m!}\sum\limits_{1\leq j_1\neq\ldots\neq j_k\leq m}\int\limits_{1}^{\h}\int\limits_{0}^{\pi}\ldots\int\limits_{0}^{\pi}J_f(y,\theta_1,\ldots,\theta_{m})\prod\limits_{i=1}^{k}\ind_{I_i}(\theta_{j_i})\dd y\,\dd\theta_1\ldots\dd\theta_{m}.
\end{align*}
Using again the invariance with respect to permutations of $\theta_1,\ldots,\theta_{m}$ we conclude
\[
W=\frac{1}{(m-k)!}\int\limits_{1}^{\h}\int\limits_{0}^{\pi}\ldots\int\limits_{0}^{\pi}J_f(y,\theta_1,\ldots,\theta_{m})\prod\limits_{i=1}^{k}\ind_{I_i}(\theta_{i})\dd y\,\dd\theta_1\ldots\dd\theta_{m}.
\]

Applying Lemma \ref{lm2} we finally obtain
\begin{equation}\label{eq_24}
W=\kappa_m\,\h^{m+1}+O(\h^{m}),
\end{equation}
where
\begin{align}
\kappa_m:&=\frac{2^{m(m+1)/2}}{(m+1)(m-k)!}\int\limits_{0}^{\pi}\ldots\int\limits_{0}^{\pi}\prod\limits_{1\leq i<j\leq m}\left|\theta_i-\theta_j\right|\prod\limits_{l=1}^{m}\sin\theta_l\prod\limits_{1}^{k}\ind_{I_i}(\theta_{i})\dd \theta_1\ldots\dd \theta_{m}\notag\\
&=\frac{2^{m(m+1)/2}}{(m+1)(m-k)!}\int\limits_{I_1}\ldots\int\limits_{I_k}\int\limits_{0}^{\pi}\ldots\int\limits_{0}^{\pi}\prod\limits_{1\leq i<j\leq m}\left|\theta_i-\theta_j\right|\prod\limits_{l=1}^{m}\sin\theta_l\dd \theta_1\ldots\dd \theta_{m}\notag\\
&=\omega_m\int\limits_{I_1}\ldots\int\limits_{I_k}\rho_{m,k}(\theta_1,\ldots,\theta_k)\dd \theta_1\ldots\dd \theta_{k}, \label{eq_21}
\end{align}
with
\begin{equation}\label{eq_22}
\rho_{m,k}(\theta_1,\ldots,\theta_k):=Z_{m}^{-1}\frac{m!}{(m-k)!}\int\limits_{0}^{\pi}\ldots\int\limits_{0}^{\pi}\prod\limits_{l=1}^{m}\sin\theta_l\prod\limits_{1\leq i<j\leq m}\left|\theta_i-\theta_j\right|\dd \theta_{k+1}\ldots\dd \theta_{m}.
\end{equation}
and
\[
Z_{m}=2^{m(m+1)/2}S_{m}\left(1,1,1/2\right).
\]
 
\begin{lemma}\label{lm5}
For any integer $m$ and $1\leq k\leq m$ we have 
	\[
	\rho_{m,k}(\theta_1,\ldots,\theta_{k})=\prod\limits_{l=1}^{k}\sin\theta_l\Pf\Big[K_m(\cos\theta_i,\cos\theta_j)\Big]_{i,j=1,\ldots,k},
	\]
	where $K_m(x,y)$ is defined by \eqref{eq_23} - \eqref{eq_58}.
\end{lemma}

The proof is given in Subsection \ref{proof_lm5}.

Combining Lemma \ref{lm5} with equations \eqref{eq_24}, \eqref{eq_21} and \eqref{eq_18} we finish the proof of Theorem \ref{th2}.

\section{Proofs of Lemmas}

\subsection{Proof of Lemma \ref{lm3}}\label{proof_lm3}

Recall that $V_m$ denotes the set of points $\ba\in\R^{m+1}$ such that the roots of $P_{\ba}(t)$ may be written as \eqref{eq_1}. Furthermore, recall that $V^l_{m,k}\subset V^m$ denotes the subset of points $\ba\in V^m$ such that there are exactly $l$ tuples $(e^{i\theta_1},\ldots,e^{i\theta_k})$ such that  
\[
\theta_j\in I_j:=[s_{2j-1};s_{2j}],
\]
and $e^{i\theta_j}$, $1\leq j\leq k$ are the roots of $P_{\ba}(t)$.
It is clear that for fixed $k$ and $m$ all sets $V^l_{m,k}$ are disjoint, only finite number, say $L$, of them are non-empty and, moreover,
\[
V_m=\bigcup\limits_{l=0}^{L}V^l_{m,k}.
\]
Thus, the boundary of the set $V^l_{m,k}$ can be covered by $\partial V_m$ and a set
\[
H:=\left\{\ba\in\R^{m+1}\colon P_{\ba}\left(e^{is_j}\right)=0 \text{ for at least one } 1\leq j\leq 2k\right\}.
\]

First of all recall that according to Lemma \ref{lm1} there exists a bijective map $f: (1,\h]\times L_{m}\rightarrow V_m$ defined by \eqref{eq_3} and \eqref{eq_2}, where the simplex $L_{m}$ is defined by \eqref{eq_26}.
This map defines a homeomorphism of manifolds with boundary, which means that $\partial V_m = f\left(\partial \left((1,\h]\times L_m\right)\right)$, where
\[
\partial \left((1,\h]\times L_{m}\right)\subset J_1^{+}\cup J_1^{-}\cup J_2^{+}\cup J_2^{-}\cup J_3,
\]
and
\begin{align*}
J_1^{+}:&=\left\{\bx\in\R^{m+1}\colon x_0=\h \right\}, \quad J_1^{-}:=\left\{\bx\in\R^{m+1}\colon x_0=1 \right\},\\
J_2^{+}:&=\left\{\bx\in\R^{m+1}\colon x_m=\pi \right\}, \quad J_2^{-}:=\left\{\bx\in\R^{m+1}\colon x_1=0 \right\},\\
J_3:&=\left\{\bx\in\R^{m+1}\colon x_i=x_{i+1} \text{ for at least one } 1\leq i\leq m\right\}.
\end{align*}
Moreover, according to the definition for any $\ba\in H$ at least one of the equalities
\[
a_{m+1}=-2\cos \left((m+1)s_{j}\right)-2\sum\limits_{i=1}^{m}a_{m+1-i}\cos\left(is_j\right),\quad 1\leq j\leq 2k,
\]
hold.

Let us construct maps $\phi_1, \dots, \phi_{3+m+2k} : [0;1]^{m} \to\R^{m+1}$ as follows
\begin{align*}
\phi_1(t_1,\ldots,t_{m}):&= f\left(\h, \pi\,t_1, \ldots, \pi\,t_m\right),\\ 
\phi_2(t_1,\ldots,t_{m}):&= f\left(1, \pi\,t_1, \ldots, \pi\,t_m\right),\\
\phi_3(t_1,\ldots,t_{m}):&= f\left(\h\, t_1+1, 0, \pi\,t_2, \ldots, \pi\,t_m\right),\\ 
\phi_4(t_1,\ldots,t_{m}):&= f\left(\h\, t_1+1, \pi\,t_2, \ldots, \pi\,t_m, \pi\right),\\
\phi_{3+j}(t_1,\ldots,t_{m}):&= f\left(\h\, t_1+1, \pi\,t_2, \ldots, \pi\, t_j, \pi\, t_j, \ldots, \pi\,t_m\right), \quad 2\leq j\leq m,
\end{align*}
and for $1\leq j\leq 2k$
\begin{align*}
\phi_{3+m+j,i}(t_1,\ldots,t_{m}):&= wt_i,\quad 1\leq i\leq m,\\ 
\phi_{3+m+j,m+1}(t_1,\ldots,t_{m}):&= -2\cos \left((m+1)s_{j}\right)-2w\sum\limits_{i=1}^{m}t_{m+1-i}\cos\left(is_j\right),\\
w:&=\max\limits_{\ba\in V_m}\max\limits_{1\leq i\leq m+1}|a_i|.
\end{align*}

It is easy to see that all maps $\phi_i$ are Lipschitz continuous since they are continuously differentiable in a compact set. Moreover, from the arguments above it follows that
\[
\partial V^l_{m,k}\subset \bigcup\limits_{j=1}^{m+3+2k}\phi_j\left([0;1]^{m}\right).
\]

Due to the definition of mapping $f: (1,\h]\times L_{m}\rightarrow V_m$ we have $w\leq \tilde c\,\h$ for some constant $\tilde c$ depending on $m$ only and, thus, for the Lipschitz constant $L_j$ of the maps $\phi_j$ we conclude
\[
L_j:=\sup\limits_{\bx\in[0;1]^{m}}\max\limits_{1\leq i\leq m}\sum\limits_{q=1}^{m+1}\frac{\partial \phi_{j,q}}{\partial t_i}(\bx)\leq c\,\h,
\]
for some constant $c$ depending on $m$ only.

Finally we use a lattice point counting result by Widmer \cite{W10}.

\begin{theorem}
Let $\Lambda$ be a lattice in $\R^d$ with successive minima $\lambda_1,\ldots, \lambda_d$. Let $B\subset\R^d$
be a bounded set with boundary of Lipschitz class $(M,L)$. Then $B$
is measurable, and, moreover,
\[
\left|\#(\Lambda\cap B)-\frac{\Vol(B)}{\det \Lambda}\right|\leq c(d)\,M\,\max\limits_{0\leq i<d}\frac{L^i}{\lambda_1\cdots\lambda_i}.
\]
For $i = 0$ the expression in the maximum is to be understood as 1.
\end{theorem}

Taking $\Lambda = \Z^{m+1}$, $B=V_{m,k}^l$, $L=c\,\h$, and $M=m+3+2k$, and applying the theorem above we get
\[
\left|\#(\Z^{m+1}\cap V_{m,k}^l)-v_{m,k}^l\right|\leq c(m)c_m^{m}(m+3+2k)((m+1)!)^{-1}\h^{m}=:C_1\h^{m},
\]
which finishes the proof.

\subsection{Proof of Lemma \ref{lm4}}\label{proof_lm4}

Consider a reducible polynomial $P\in\Pl_m(\h)$ and assume that it can be written as a product of polynomials $P_1,P_2\in\Z[t]$,  such that $\deg P_1 = 2m_1$, $\deg P_2 = 2m_2$ and $m_1+m_2 = m+1$. By definition of $\Pl_m(\h)$ the polynomial $P$ is monic and has roots described by \eqref{eq_1}. Hence the polynomials $P_1$ and $P_2$ are monic as well, and all roots of one of them (say, $P_1$) are lying on the unit circle $\T$ and the one (say, $P_2$) belongs to the set $\Pl_{m_2-1}(\h)$. Moreover, by Kronecker's theorem \cite{Kr} we conclude that the polynomial $P_1$ has to be a product of cyclotomic polynomials.

From the arguments above it follows that $\#\Pl^{R}_m(\h)$ does not exceed the number of pairs $\left(P_1,P_2\right)\in\Z[t]\times\Z[t]$ of monic polynomials with integer coefficients such that $P_1$, $\deg P_1=2m_1$ is a product of cyclotomic polynomials, $P_2\in\Pl_{m_2-1}(\h)$ and $m_1+m_2=m+1$. It is easy to see that the number of polynomials $P_1$ does not exceed some constant $c_1(m_1)$, and from Lemma \ref{lm3} and \eqref{eq_11} it immediately follows that
\[
\#\Pl_m(\h)=\frac{2^{m(m+1)}}{m+1}\,\prod\limits_{k=0}^{m-1}\frac{k!^2}{(2k+1)!}\,\h^{m+1}+O(\h^{m})<c_2(m)\,\h^{m+1}.
\]
Hence, we conclude
\[
\#\Pl^{R}_m(\h)\leq \sum\limits_{k=1}^{m}c_1(k)\,\#\Pl_{m-k}(\h)\leq \sum\limits_{k=1}^{m}c_1(k)c_2(m-k)\h^{m+1-k}\leq C_2\h^{m},
\]
where the constant $C_2$ depends on $m$ only.

\subsection{Proof of Lemma \ref{lm1}}\label{proof_lm1}

Let us consider a polynomial $P_{\ba}(t)$ with roots \eqref{eq_1}. Then
\begin{align*}
P_{\ba}(t)&=\left(t-y\right)\left(t-y^{-1}\right)\left(t-e^{i\theta_1}\right)\left(t-e^{-i\theta_1}\right)\cdots\left(t-e^{i\theta_{m}}\right)\left(t-e^{-i\theta_{m}}\right)\\
&=\left(t^2-(y+y^{-1})t+1\right)\left(t^2-2t\cos\theta_1+1\right)\cdots\left(t^2-2t\cos\theta_{m}+1\right)\\
&=\left(t^2-z_0t+1\right)\left(t^2-z_1t+1\right)\cdots\left(t^2-z_{m}t+1\right)\\
&=t^{2m+2}+a_1t^{2m+1}+\ldots+a_{m+1}t^{m+1}+\ldots+a_1t+1.
\end{align*}
All possible summands in the final expression consist of products of one of the summands of $t^2$, $1$ or $-z_kt$ from every parentheses.

Consider the summand $a_{2i} t^{2i}$. The power $t^{2i}$ appears in the following combinations only 
\[
\underbrace{t^2,\ldots, t^2}_\text{$i-j$ times},-z_{k_1}t,\ldots, -z_{k_{2j}}t, \underbrace{1,\ldots, 1}_\text{$m+1-i-j$ times},
\]
where $0\leq j\leq i$ and $0\leq k_1<\ldots<k_{2j}\leq m$.
Hence it follows that $a_{2i}$ is the sum of the following terms
\[
{m+1-2j \choose i-j }\sum\limits_{0\leq k_1<\ldots <k_{2j}\leq m}\prod_{s=1}^{2j}z_{k_{s}},
\]
which leads to the formula \eqref{eq_2}.

In the same manner consider the summand $a_{2i-1} t^{2i-1}$ and note that the power $t^{2i-1}$ appears in the following combinations only 
\[
\underbrace{t^2,\ldots, t^2}_\text{$i-j-1$ times},-z_{k_1}t,\ldots, -z_{k_{2j+1}}t, \underbrace{1,\ldots, 1}_\text{$m+1-i-j$ times},
\]
where $0\leq j\leq i-1$ and $0\leq k_1<\ldots<k_{2j+1}\leq m$. Thus, $a_{2i-1}$ is the sum of the following terms
\[
-{m-2j \choose i-j-1 }\sum\limits_{0\leq k_1<\ldots <k_{2j+1}\leq m}\prod_{s=1}^{2j+1}z_{k_{s}},
\]
which leads to the formula \eqref{eq_3}.

\subsection{Proof of Lemma \ref{lm2}}\label{proof_lm2}

By definition the Jacobian $J_f(y,\theta_1,\ldots,\theta_{m})$ is the absolute value of the determinant of the following matrix
\[
A:=\left(
\begin{matrix}
\frac{\partial a_1}{\partial z_0}\frac{\partial z_0}{\partial y} & \frac{\partial a_1}{\partial z_1}\frac{\partial z_1}{\partial \theta_1} & \dots & \frac{\partial a_1}{\partial z_{m}}\frac{\partial z_{m}}{\partial \theta_{m}}\\
\frac{\partial a_2}{\partial z_0}\frac{\partial z_0}{\partial y} & \frac{\partial a_2}{\partial z_1}\frac{\partial z_1}{\partial \theta_1} & \dots & \frac{\partial a_2}{\partial z_{m}}\frac{\partial z_{m}}{\partial \theta_{m}}\\
\vdots & \vdots & \ddots & \vdots \\
\frac{\partial a_{m}}{\partial z_0}\frac{\partial z_0}{\partial y} & \frac{\partial a_{m}}{\partial z_1}\frac{\partial z_1}{\partial \theta_1} & \dots & \frac{\partial a_{m}}{\partial z_{m}}\frac{\partial z_{m}}{\partial \theta_{m}}
\end{matrix}
\right).
\]
Hence,
\begin{equation}\label{eq_7}
\left|\det(A)\right|=\left|\frac{\partial z_0}{\partial y}\,\prod\limits_{l=1}^{m}\frac{\partial z_l}{\partial \theta_l}\right|\left|\det(\tilde A)\right|=2^{m}\left|1-y^{-2}\right|\prod\limits_{l=1}^{m}\left|\sin\theta_l\right|\cdot\Delta,
\end{equation}
where 
\[
\Delta=\det\left(
\begin{matrix}
1 & \dots & 1\\
\sum\limits_{k_1\neq 0}z_{k_1} & \dots & \sum\limits_{k_1\neq m-1}z_{k_1} \\
\vdots& \ddots & \vdots \\
\sum\limits_{j=0}^{i-1}{m-2j \choose i-j-1 }\sum\limits_{k_1, \ldots, k_{2j}\neq 0}\prod\limits_{s=1}^{2j}z_{k_{s}} & \dots & \sum\limits_{j=0}^{i-1}{m-2j \choose i-j-1 }\sum\limits_{k_1, \ldots, k_{2j}\neq m}\prod\limits_{s=1}^{2j}z_{k_{s}}\\
\sum\limits_{j=1}^{i}{m+1-2j \choose i-j}\sum\limits_{k_1,\ldots,k_{2j-1}\neq 0}\prod\limits_{s=1}^{2j-1}z_{k_{s}} & \dots & \sum\limits_{j=1}^{i}{m+1-2j \choose i-j}\sum\limits_{k_1,\ldots,k_{2j-1}\neq m}\prod\limits_{s=1}^{2j-1}z_{k_{s}}\\
\vdots
\end{matrix}
\right).
\]
We would like to show that 
\begin{equation}\label{eq_6}
\Delta=\det\left(
\begin{matrix}
1 & \dots & 1\\
\sum\limits_{k_1\neq 0}z_{k_1} & \dots & \sum\limits_{k_1\neq m}z_{k_1} \\
\sum\limits_{k_1,k_2\neq 0}z_{k_1} z_{k_2}& \dots & \sum\limits_{k_1,k_2\neq m}z_{k_1}z_{k_2} \\
\vdots& \ddots & \vdots \\
\prod\limits_{s\neq 0}z_{s} & \dots & \prod\limits_{s\neq m}z_{s}
\end{matrix}
\right).
\end{equation}
We shall prove this formula by induction. First off all notice that the first and the second rows already have the required form, hence we may take them as base of the induction. Assume all rows up to $(p-1)$-th row  are of the form in \eqref{eq_6}. Consider the $p$-th row, which has the form $(a_{p,0},\ldots, a_{p,m})$, where
\[
a_{p,l}=\sum_{j=0}^{p-2}c_{j,p}\sum\limits_{k_1,\ldots,k_{j}\neq l}\prod\limits_{s=1}^{j}z_{k_{s}}+\sum\limits_{k_1,\ldots,k_{p-1}\neq l}\prod\limits_{s=1}^{p-1}z_{k_{s}},
\]
with the constants $c_{j,p}\ge 0$ depending on $j$, $p$, and $m$ only. Multiplying the $i$-th row by $c_{i-1,p}$ and subtracting it from the $p$-th row for all $1\leq i\leq p-1$ we have
\[
\tilde a_{p,l}=\sum\limits_{k_1,\ldots,k_{p-1}\neq l}\prod\limits_{s=1}^{p-1}z_{k_{s}},
\]
and, since these transformations do not change the determinant of the matrix \eqref{eq_6} holds for the $p$-th row as well.

On the other hand completing the sums in the rows and applying successive row operations leaving the determinant $\Delta$ invariant  we finally arrive at
\[
\Delta=	\det\left(
\begin{matrix}
1 & \dots & 1\\
z_{0} & \dots & z_{m} \\
\vdots& \ddots & \vdots \\
z^{m}_{0}& \dots & z^{m}_{m}
\end{matrix}\right).
\]
The last determinant is the well-known Vandermonde determinant and, thus, we have
\[
\Delta=2^{m(m-1)/2}\prod\limits_{l=1}^{m}\left|y+y^{-1}-2\cos\theta_i\right|\prod\limits_{1\leq i<j\leq m}\left|\cos\theta_i-\cos\theta_j\right|.
\]
Substituting this into \eqref{eq_7} and taking into account that $y>1$ yields Lemma \ref{lm2}.

\subsection{Proof of Lemma \ref{lm5}}\label{proof_lm5}

Consider the expression 
\begin{equation}\label{eq_40}
U:=\int\limits_{I_1}\ldots\int\limits_{I_k}\rho_{m,k}(\theta_1,\ldots,\theta_k)\dd \theta_1\ldots\dd \theta_{k}.
\end{equation}

Using change of variables $x_i=-\cos\theta_i$, $1\leq i\leq m$ and the definition \eqref{eq_22} of the function $\rho_{m,k}(\theta_1,\ldots,\theta_k)$ we write
\[
U=Z_{m}^{-1}\frac{m!}{(m-k)!}\int\limits_{T_1}\ldots\int\limits_{T_k}\int\limits_{-1}^{1}\ldots\int\limits_{-1}^{1}\prod\limits_{1\leq i<j\leq m}\left|x_i-x_j\right|\dd x_{1}\ldots\dd x_{m},
\]
where 
\[
T_i:=\left\{-\cos\theta\colon \theta\in I_i\right\}\subset [-1;1],\quad 1\leq i\leq k.
\]

Define the following function
\[
p_m(x_1,\ldots,x_m):=\prod\limits_{1\leq i<j\leq m}\left|x_i-x_j\right|.
\]
According to \cite[Chapter 19]{Mehta} this function corresponds to the joint probability density function of the eigenvalues of the Jacobi random matrix ensemble with $\beta=1$ and weight function $w(x)=1$, $-1\leq x\leq 1$ (see Appendix \ref{theory} for the more details).

Then
\[
R_k(x_1,\ldots,x_k):=Z_{m}^{-1}\frac{m!}{(m-k)!}\int\limits_{-1}^{1}\ldots\int\limits_{-1}^{1}p_m(x_1,\ldots,x_m)\dd x_{k+1}\ldots\dd x_{m}
\]
defines a $k$-point correlation function (see \eqref{eq_50}), which in this case can be written in the following form
\[
R_k(x_1,\ldots,x_k)=\Pf[K_N(x_i,x_j)]_{i,j=1,\ldots,k}
\]
with the Kernel function $K_N(x,y)$ is defined by \eqref{eq_23} - \eqref{eq_58} (see Apendix \ref{jacobi}).

Finally we write
\begin{multline*}
U=\int\limits_{T_1}\ldots\int\limits_{T_k}R_k(x_1,\ldots,x_k)\dd x_1\ldots\dd x_k=\int\limits_{T_1}\ldots\int\limits_{T_k}\Pf[K_N(x_i,x_j)]_{i,j=1,\ldots,k}\dd x_1\ldots\dd x_{k}\\
=\int\limits_{I_1}\ldots\int\limits_{I_k}\prod\limits_{l=1}^{k}\sin\theta_l\Pf\Big[K_N(-\cos\theta_i,-\cos\theta_j)\Big]_{i,j=1,\ldots,k}\dd \theta_1\ldots\dd \theta_{k},
\end{multline*}
and combining this with \eqref{eq_40} we finish the proof.

\section*{Acknowledgments}

The authors are grateful to Arturas Dubickas for the suggestion of the problem and to Zakhar Kabluchko and Gernot Akemann for helpful discussions which improved the result.

\appendix

\section{Random Matrices: General Theory}\label{theory}

In this appendix we collect some basic facts about random matrices and the distribution of their eigenvalues. In view of the rather extensive literature on this subject we shall restrict ourselves to a very short review of results we need.

One way of defining  a random matrix ensemble is to specify the joint probability density functions for its eigenvalues in the following form
\[
p^{\beta}_N(x_1,\ldots,x_N)=\prod\limits_{i=1}^N w(x_i)\prod\limits_{1\leq i<j\leq N}\left|x_i-x_j\right|^{\beta},
\]
where $\beta$ is an in general complex parameter and the so called weight function $w(x)$ can be chosen to suit the needs. The most well-studied cases are $\beta =1$ (real ensembles), $\beta=2$ (complex ensembles) and $\beta=4$ (quaternion ensembles).

One of the main objectives of Random Matrix Theory is to investigate the distribution of eigenvalues of different random matrix ensembles and particularly their limiting distribution when $N\rightarrow \infty$. For this purpose one needs to calculate $k$-point correlation functions of eigenvalues \cite[eq. (5.7.1)]{Mehta} defined by
\begin{equation}\label{eq_50}
R_k^{\beta}(x_1,\ldots,x_k):=Z_{N,\beta}^{-1}\frac{N!}{(N-k)!}\int p^{\beta}_N(x_1,\ldots,x_N)\dd x_{k+1}\ldots\dd x_N,
\end{equation}
where the normalization constant $Z_{N,\beta}$ is given by
\[
Z_{N,\beta}=\int p^{\beta}_N(x_1,\ldots,x_N)\dd x_{1}\ldots\dd x_N.
\]

The eigenvalues of random matrix ensemble form a point process (see \cite{HKPV09} for a definition and theory of point processes) and for two special types of point processes, namely the  {\it Determinantal point processes} and the {\it Pfaffian point processes} the  functions \eqref{eq_50} will have a compact representation.
\begin{enumerate}
	\item In case of Determinantal point processes all $k$-point correlation functions have the form
	\[
	R_k^{\beta}(x_1,\ldots,x_k)=\det\left[K_N(x_i,x_j)\right]_{i,j=1,\ldots,k};
	\]
	\item In case of Pfaffian point processes all $k$-point correlation functions have the form
	\[
	R_k^{\beta}(x_1,\ldots,x_k)=\Pf\left[K_N(x_i,x_j)\right]_{i,j=1,\ldots,k},
	\]
\end{enumerate}
where $K_N(x,y)$ is some special function satisfying the conditions described by \cite[eq. (5.1.21), (5.1.22)]{Mehta}. The function $K_N(x,y)$ is called a {\it Kernel function}.

It is known from random matrix theory that the case $\beta=2$ corresponds to Determinantal point processes and the case  $\beta=1$, $\beta=4$ correspond to Pfaffian point processes. Formulas for the kernel function may be derived using systems of orthogonal ($\beta=2$) and skew-orthogonal ($\beta=1$, $\beta=4$) polynomials with the particular choice depending on the weight function $w(x)$. We skip the details referring the reader to \cite[Chapter 5]{Mehta}.

\section{Jacobi Ensemble with $\beta=1$}\label{jacobi}

This appendix is devoted to the particular random matrix ensemble used in the proof of Theorem \ref{th2}, namely the Jacobi $\beta$-ensemble \cite[Chapter 19]{Mehta} with $\beta=1$.

Consider the random matrix ensemble with the following joint probability density function for eigenvalues
\[
p_N(x_1,\ldots,x_N)=\prod\limits_{1\leq i<j\leq N}\left|x_i-x_j\right|,\quad x_i\in[-1;1],
\]
with the weight function
\begin{equation}\label{eq_52}
w(x)=\begin{cases}
1,\quad -1\leq x\leq 1,\\
0, \quad\text{otherwise}.
\end{cases}
\end{equation}
This is a special case of more general class of weight functions
\[
w(x):=w(a,b;x)
=\begin{cases}
(1-x)^a(1+x)^b,\quad -1\leq x\leq 1,\\
0, \quad\text{otherwise},
\end{cases}
\quad a,b>-1,
\]
which gives an ensemble corresponding to the Jacobi orthogonal polynomials $\jac{j}{a}{b}(t)$ defined by \eqref{eq_51}. The orthogonality property means
\[
\int\limits_{-1}^{1}w(a,b;x)\jac{j}{a}{b}(x)\jac{i}{a}{b}(x)\dd x=h_j^{(a,b)}\delta_{i,j},
\]
where
\[
h_j^{(a,b)}=\frac{2^{a+b+1}}{(2j+a+b+1)}\frac{\Gamma(j+a+1)\Gamma(j+b+1)}{\Gamma(j+1)\Gamma(j+a+b+1)}.
\]
The Jacobi polynomials have been known for a long time and their properties are well studied \cite{Szego}.

As it was mentioned in Appendix \ref{theory} for $\beta=1$ one needs to construct the set $R_j(t)$ of skew-orthogonal polynomials corresponding to the weight function \eqref{eq_52}. Skew orthogonality 
means that
\begin{align*}
\int\limits_{-1}^{1}R_{2j}(x)R_{2i+1}(x)\dd x&=r_j\delta_{i,j},\\
\int\limits_{-1}^{1}R_{2j}(x)R_{2i}(x)\dd x&=\int\limits_{-1}^{1}R_{2j+1}(x)R_{2i+1}(x)\dd x=0,
\end{align*}

Thus, following the procedure described in \cite[Chapter 19.2]{Mehta} and using the identities for Jacobi polynomials
\begin{align*}
\frac{\dd}{\dd\, t}\jac{j}{a}{b}(t)&=\frac{\Gamma(a+b+j+2)}{2\Gamma(a+b+j+1)}\jac{j-1}{a+1}{b+1}(t),\\
\jac{j}{a}{b}(-z)&=(-1)^j\jac{j}{a}{b}(z),
\end{align*}
we get for $N=2s$
\begin{align}
R_{2j}(t)&=\jac{2j}{1}{1}(t),\notag\\
R_{2j+1}(t)&=\frac{\dd}{\dd\,t}\left((t^2-1)\jac{2j}{1}{1}(t)\right),\notag\\
r_j&=\frac{8(2j+1)}{(4j+3)(2j+2)},\label{eq_59}
\end{align}
and for $N=2s+1$
\begin{align}
R_{2j}(t)&=\jac{2j+1}{1}{1}(t),\notag\\
R_{2j+1}(t)&=\frac{\dd}{\dd\,t}\left((t^2-1)\jac{2j+1}{1}{1}(t)\right),\notag\\
R_{2s}(t)&=\frac{s+1}{2}\jac{2s}{1}{1}(t),\notag\\
r_j&=\frac{8(2j+2)}{(4j+5)(2j+3)}.\label{eq_60}
\end{align}

Define the function
\[
\psi_j(t):=\frac12\int\limits_{-1}^{1}\sign(t-x)R_j(x)\dd\,x.
\]
Thus for $N=2s$ we have
\begin{align}
\psi_{2j}(t)&=\frac{1}{j+1}\jac{2j+1}{0}{0}(t),\label{eq_61}\\
\psi_{2j+1}(t)&=(t^2-1)\jac{2j}{1}{1}(t),\label{eq_62}
\end{align}
and for $N=2s+1$ we get
\begin{align}
\psi_{2j}(t)&=\frac{2}{2j+3}\left(\jac{2j+2}{0}{0}(t)-1\right),\label{eq_63}\\
\psi_{2j+1}(t)&=(t^2-1)\jac{2j+1}{1}{1}(t).\label{eq_64}
\end{align}
Finally the Kernel function $K_N(x,y)$ for the Jacobi ensemble with $\beta=1$ is defined via equations (19.2.22) - (19.2.28) in \cite{Mehta}.

Note that in deriving \eqref{eq_53} - \eqref{eq_55} we  have combined the equations (19.2.23) - (19.2.26) with (19.2.27) and (19.2.28). Moreover using the notation $c:=\left(N\mod{2}\right)$ and equations \eqref{eq_59} - \eqref{eq_64} we arrive at \eqref{eq_56} - \eqref{eq_58}. Furthermore, note that the representation of $K_N(x,y)$ in (19.2.22) is given in quaternion form. In order to get a representation for Pfaffian point processes we have to use \cite[Theorem 5.1.2]{Mehta}.

\bibliographystyle{plainnat}

\begin{thebibliography}{10}
	
	\bibitem{W10}
	M.~Widmer
	\newblock Counting primitive points of bounded height.
	\newblock {\em Trans. Amer. Math. Soc.}, 362:4793--4829, 2010.
	
	\bibitem{CH17}
	F.~Calegari and Z.~Huang
	\newblock Counting Perron numbers by absolute value.
	\newblock {\em J. London Math. Soc.}, (96):181--200, 2017.
	
	\bibitem{hD51}
	H.~Davenport.
	\newblock On a principle of {L}ipschitz.
	\newblock {\em J.~Lond. Math. Soc.}, 26(3):179--183, 1951.
	\newblock Corrigendum: ``On a principle of Lipschitz'', {\it J.~Lond. Math.
		Soc.} {\bf 39} (1964), 580.
	
	\bibitem{MV08}
	D.~Masser and J.~D. Vaaler.
	\newblock Counting algebraic numbers with large height~{I}.
	\newblock In {\em Diophantine Approximation}, volume~16 of {\em Dev. Math.},
	pages 237--243. SpringerWienNewYork, Vienna, 2008.
	
	\bibitem{Ku09}
	G. Kuba.
	\newblock On the distribution of reducible polynomials.
	\newblock {\em Math. Slovaca.}, 59(3):349--356, 2009.
	
	\bibitem{Prasolov}
	V.V. Prasolov.
	\newblock Polynomials.
	\newblock vol. 11 of Algorithms and Computation in Mathematics, Springer, Berlin, 2004.
	
	\bibitem{mW12}
	M.~Widmer.
	\newblock Lipschitz class, narrow class, and counting lattice points.
	\newblock {\em Proc. Amer. Math. Soc.}, 140(2):677--689, 2012.
	
	\bibitem{rL65}
	R.~Lipschitz.
	\newblock {\"U}ber die asymptotischen {G}esetze von gewissen {G}attungen
	zahlentheoretischer {F}unktionen.
	\newblock {\em Monatsber. der Berliner Akademie}, 1865:174--185, 1865.
	
	\bibitem{Kr}
	L. Kronecker.
	\newblock Zwei S\"{a}tze \"{u}ber Gleichungen mit ganzzahligen Coefficienten.
	\newblock {\em J. reine angew. Math.}, 53: 173---175, 1857.
	
	\bibitem{Leh}
	D. H. Lehmer.
	\newblock Factorization of certain cyclotomic functions.
	\newblock {\em Ann. of Math.}, 34: 461---469, 1933.
	
	\bibitem{PiSa}
	M.-J.  Bertin,  A.  Decomps-Guilloux,  M.  Grandet-Hugot,  M.  Pathiaux-Delefosse,  and  J.P. Schreiber.
	\newblock Pisot and Salem numbers.
	\newblock Birkhauser Verlag, Basel, 1992.  
	
	\bibitem{Boy81}
	D. W. Boyd. 
	\newblock Speculations  concerning  the  range  of  Mahler's  measure.
	\newblock {\em Canad.  Math.  Bull.}, 24(4): 453---469, 1981.
	
	\bibitem{Sal63}
	R. Salem.
	\newblock Algebraic numbers and Fourier analysis.
	\newblock Heath, Boston, MA, 1963.

	\bibitem{EvWar}
	G. Everest and T. Ward.
	\newblock Heights of polynomials and entropy in algebraic dynamics.
	\newblock Universitext, Springer-Verlag London, Ltd., London, 1999. 
	
	\bibitem{Sel44}
	A. Selberg.
	\newblock Remarks on a multiple integral.
	\newblock {\em Norsk Mat. Tidsskr.}, 26: 71--78, 1944.
	
	\bibitem{Szego}
	G. Szeg\"o.
	\newblock Orthogonal Polynomials. 
	\newblock Colloquium Publications. XXIII. American Mathematical Society, 1939.
	
	\bibitem{Mehta}
	M.L. Mehta.
	\newblock Random Matrices. 
	\newblock Pure and Applied Mathematics (Vol 142). Elsevier, 2004.
	
	%
	\bibitem{HKPV09}
	J.B. Hough, M. Krishnapur, Y. Peres, and B. Virag. 
	\newblock Zeros of Gaussian Analytic Functions and Determinantal Point Processes.
	\newblock American Mathematical Society, Providence (RI), 2009.
	
	\bibitem{Bar14}
	F. Barroero.

	\newblock Counting algebraic integers of fixed degree and bounded height.
	\newblock  {\em Monatsh. Math.},
175(1):25--41, 2014.
	
	\bibitem{GG17}
	R. Grizzard and J. Gunther.
	\newblock Slicing the stars: counting algebraic numbers, integers, and units by degree and height.
	\newblock {\em Algebra Number Theory}, 11(6): 1385--1436, 2017.

    \bibitem{vB99}
	V. Beresnevich. 
	\newblock On approximation of real numbers by real algebraic numbers.
	\newblock {\em Acta Arithmetica}, 90(2):97--112, 1999.

	\bibitem{Sch93}
	W.M. Schmidt.
	\newblock Northcott's theorem on heights I. A general estimate.,
	\newblock {\em Monatsh. Math.}, 115: 169--181, 1993.

	\bibitem{MV08}
	D.~Masser and J.~D. Vaaler.
	\newblock Counting algebraic numbers with large height~{I}.
	\newblock In {\em Diophantine Approximation}, volume~16 of {\em Dev. Math.},
	pages 237--243. SpringerWienNewYork, Vienna, 2008.
	
	\bibitem{dK14}
	D.~Kaliada.
	\newblock On the density function of the distribution of real algebraic
	numbers.
	\newblock {\em Journal de Th\'eorie des Nombres de Bordeaux}, 29(1): 179--200, 2017.
	
	\bibitem{GKZ15}
	F.~G{\"o}tze, D.~Kaliada, and D.~Zaporozhets.
	\newblock Distribution of complex algebraic numbers.
	\newblock {\em Proc. Amer. Math. Soc.}, 145(1):61--71, 2017.
	\newblock Preprint arXiv:1410.3623, 2014.
	
	\bibitem{GKZ17}
	F.~G{\"o}tze, D.~Kaliada, and D.~Zaporozhets.
	\newblock Joint distribution of conjugate algebraic numbers: a random polynomial approach.
	\newblock Preprint arXiv:1703.02289, 2017.

\end{thebibliography}

\end{document}